\documentclass[11pt]{article}

\usepackage{amsmath,amssymb,amscd,amsthm}
\usepackage[all]{xy}
\usepackage{enumerate}
\usepackage[dvips]{graphicx}
\usepackage{wrapfig}
\theoremstyle{plain}
\numberwithin{equation}{section}
\newtheorem{thm}{Theorem}[section]
\newtheorem{prop}[thm]{Proposition}
\newtheorem{cor}[thm]{Corollary}
\newtheorem{lem}[thm]{Lemma}

\theoremstyle{definition}
\newtheorem{dfn}[thm]{Definition}

\newtheorem{rmk}[thm]{Remark}

\def\rank{\mathop{\mathrm{rank}}\nolimits}
\def\dim{\mathop{\mathrm{dim}}\nolimits}

\def\Hom{\mathop{\mathrm{Hom}}\nolimits}

\def\<{{\langle}}
\def\>{{\rangle}}

\def\Aut{\mathop{\mathrm{Aut}}\nolimits}

\newcommand{\bb}[1]{{\mathbb{#1}}}
\newcommand{\mca}[1]{{\mathcal{#1}}}
\newcommand{\mr}[1]{{\mathrm{#1}}}

\title{The group of autoequivalences and the Fourier-Mukai number of a projective manifold}
\author{Kotaro Kawatani}
\date{Osaka University}

\begin{document}
\maketitle

\begin{abstract}
Let $X$ be a smooth projective variety and $\Aut (D(X))$ the group of autoequivalences of the derived category of $X$. 
In this paper we show that $X$ has no Fourier-Mukai partner other than $X$ when $\Aut (D(X))$ is generated by shifts, automorphisms and tensor products of line bundles. 
\end{abstract}

\section{Introduction}

In this paper, a \textit{projective manifold} means a smooth projective variety  over the complex number field $\bb C$. We consider the derived category $D(X)$ of a projective manifold $X$. 
That is $D(X)$ is the bounded derived category of the abelian category $\mr{Coh}(X)$ of coherent sheaves on $X$. 
As is well-known since \cite{Muk}, for some $X$, there is another projective manifold $Y$ which is not isomorphic to $X$ but the derived category $D(Y)$ is equivalent to $D(X)$ as triangulated categories. We call such a $Y$ a \textit{Fourier-Mukai partner} of $X$. 

Let $\mr{FM}(X)$ be the set of isomorphic classes of Fourier-Mukai partners of $X$.  
It is conjectured in \cite{Kaw} that the set $\mr{FM}(X)$ is a finite set. 
For instance if $X$ is an algebraic surface, the conjecture holds (cf.\ \cite{B-M}). 
Hence we call the cardinality of $\mr{FM}(X)$ the \textit{Fourier-Mukai number} of $X$.
We note that, when $X$ is a projective K3 surface, \cite{HLOY} makes the counting formula of the Fourier-Mukai number of $X$. 

Let $\Aut (D(X))$ be the group of autoequivalences of $D(X)$. 
There are a few cases when the group $\Aut (D(X))$ is exactly determined.  
For instance $\Aut (D(X))$ is determined by \cite[Theorem 3.1]{B-O} when the canonical bundle $K_X$ (or $-K_X$) is ample. 
When $X$ is an abelian variety, $\Aut (D(X))$ is determined by \cite{Orl2}.  

In this paper we shall strudy the relation between $\Aut (D(X))$ and $\rm{FM}(X)$:
 \begin{center}
 Does $\Aut (D(X))$  give us informations on $\mr{FM}(X)$ ?
\end{center}

We shall give an answer to this question in an easy case. 
Namely our theorem is the following:

\begin{thm}\label{mainthm}(= Theorem \ref{3.2})
Let $X$ be a projective manifold. We assume that $\Aut (D(X))$ is trivially generated (cf.\ Definition \ref{2.1}). Then $\mr{FM}(X) = \{ X \}$. 
\end{thm}

As an application of Theorem \ref{mainthm}, we show the following:

\begin{cor}(= Corollary \ref{todanocor})
Let $X$ be a projective manifold such that $\deg K_X|_C\neq 0$ for any irreducible curve $C \subset X$. 
Then $\mr{FM}(X)=\{ X \}$. 
\end{cor}

We found the paper \cite{Fav} and 
noticed that Theorem \ref{mainthm} is a special case of \cite[Corollary 4.3]{Fav} on arXiv after we have finished this paper. 
However our proof is independent of Favero's proof and much simpler than his. 
In addition our motivation is essentially different from his. 

\vspace{5pt}

\noindent
\textbf{Acknowledgement}. 
I thanks {Akira Fujiki}, {Yoshinori Namikawa} and {Keiji Oguiso} who read an earlier draft and suggested several improvements. 
I'm very grateful to {Yukinobu Toda} for answering my question and giving me a  suggestion for Proposition \ref{toda}. 
\section{The group of autoequivalences}
 
In this section we recall some known results on $\Aut( D(X))$. 
First we give the following easy examples of autoequivalences:
\begin{enumerate}
\item The shift of complexes $[1]:D(X) \to D(X)$.
\item The (right derived) functor $\bb{R}f_*  = f_*: D(X) \to D(X)$, where $f$ is an automorphism of $X$.  
\item The tensor products by $L \otimes : D(X) \to D(X)$ where $L \in \mr{Pic}(X)$ and $\mr{Pic }(X)$ is the Picard group of $X$. 
\end{enumerate}

\begin{dfn}\label{2.1}
We define the subgroup $\mr{Tri}(X)$ of  $\Aut (D(X))$ by the following condition: $\mr{Tri}(X)$ is the subgroup generated by shifts, automorphisms and tensor  products with line bundles. $\mr{Tri}(X)$ is called the \textit{trivial group generated by} $X$. 
If $\Aut (D(X)) = \mr{Tri}(X)$, $\Aut (D(X))$ is said to be \textit{trivially generated}. 
\end{dfn}

\begin{rmk}
The trivial group $\mr{Tri}(X)$ is written by the follwoing:
\[
\mr{Tri}(X) = (\Aut (X) \ltimes \mr{Pic}(X)  ) \times \bb Z[1],
\]
where $\Aut (X)$ is the group of automorphisms of $X$. 
Namely for $f \in \Aut (X) $ and $L \in \mr{Pic}(X)$, we have 
\[
f_* \circ L\otimes \circ f_*^{-1} (?) = f_*(L) \otimes (?).
\]
\end{rmk}

For instance, when $K_X$ or $-K_X$ is ample, $\Aut (D(X))$ is trivially generated by \cite[Theorem 3.1]{B-O}. 
The first nontrivial example of an autoequivalence was found by Mukai. 
Let us recall his example. 

Let $A$ be an abelian variety, $\hat A$ the dual abelian variety of $A$ and $\mca P$ the Poincar\'{e} line bundle on $A \times \hat A$. 
We define the functor $\Phi:D(\hat A) \to D(A)$ by the following way:
\begin{equation}
\Phi:D(\hat A) \to D(A),\ \Phi(?) := \bb R{\pi _A}_{*}(\mca P \otimes \pi _{\hat A}^*(?)), \label{FM}
\end{equation}
where $\pi_A:A \times \hat A \to A$ and $\pi _{\hat A} : A \times \hat A \to \hat A$ are the natural projections. 
Then $\Phi$ is an equivalence between $D(\hat A)$ and $D(A)$ by \cite[Theorem 2.2]{Muk}. 
The definition (\ref{FM}) seems special, but the following theorem claims that it is sufficiently general. 

\begin{thm}(\cite[Theorem 2.18]{Orl})\label{2.2}
Let $X$ be a projective manifold and $Y$ a Fourier-Mukai partner of $X$. Then, for any equivalence $\Phi : D(X) \to D(Y)$, there is an object $\mca P^{\bullet } \in D(X \times Y)$ such that 
\[
\Phi (? ) = \bb R {\pi _Y}_{*}(\mca P^{\bullet } \stackrel{\bb L}{\otimes} \bb R\pi _X^*(?)), 
\]
where $\pi_X$ (resp.\ $\pi _Y$) is the natural projection from $X \times Y$ to $X$ (resp.\ $Y$). 
Moreover $\mca P^{\bullet}$ is unique up to isomorphism. 
\end{thm}

Thus we obtain the following useful corollary:

\begin{cor}\label{2.3}
Let $x_0$ be a closed point of $X$ and $\mca O_{x_0}$ the skyscraper sheaf of $x_0$. 
If $\Phi (\mca O_{x_0}) \simeq \mca O_{y_0}$ for a closed point $y_0 \in Y$, then there is a Zariski open subset $U$ of $X$ such that 
\[
x \in U \Rightarrow \exists y \in Y \mbox{ such that }\Phi (\mca O_x) \simeq \mca O_y. 
\]
In addition, assume that for all $x \in X$ there is a closed point $y \in Y$ such that $\Phi (\mca O_x) \simeq \mca O_y$.  
Then there is an isomorphism $f:X \to Y$ and $L \in \mr{Pic}(Y)$ such that $\Phi (?) \simeq L \otimes  ( f_*(?))$.
\end{cor}

\begin{proof}
See \cite[Corollary 5.23 and Corollary 6.14]{Huy}. 
\end{proof}

\section{Proof of Theorem \ref{mainthm}}

In this section we shall prove our main theorem. 
We first cite a key lemma of the proof essentially due to \cite{B-O}. 
We define the support $\mr{Supp}(E)$ of $E \in D(X)$ by 
\[
\mr{Supp}(E) = \bigcup _{i}\mr{Supp} (H^i(E)),
\]
where $H^i(E) $ is the $i$-th cohomology with respect to the t-structre $\mr{Coh}(X)$. 

\begin{lem}(\cite[Proposition 2.2]{B-O} or \cite[Lemma 4.5]{Huy})\label{3.1}
Let $X$ be a projective manifold and $E \in D(X)$. 
Assume that 
\[ 
\dim \mr{Supp}(E)=0 \mbox{ and } \Hom_{D(X)}(E,E[i])=
\begin{cases} 0 & (i < 0) \\ \bb C & (i=0) .\end{cases}
\]
Then $E$ is isomorphic to $\mca O_x[n]$ for some $x \in X$ and $n \in \bb Z$.  
\end{lem}

Lemma \ref{3.1} is essentially due to \cite{B-O}. The above formulation of the lemma is due to \cite{Huy}. 
\vspace{5pt}

Now let $X$ and $Y$ be projective manifolds and $\Phi :D(X) \to D(Y)$ an equivalence. 
Then we remark that $\Phi$ induces the natural group isomorphism $\Phi _*: \Aut (D(X)) \to \Aut (D(Y)) $ by the following way:
\[
\Phi_* : \Aut (D(X)) \to \Aut (D(Y)) ,\ \Phi_*(F) := \Phi \circ F \circ \Phi ^{-1}. 
\]

\begin{thm}(=Theorem \ref{mainthm})\label{3.2}
Let $X$ be a projective manifold. Assume that $\Aut (D(X))$ is trivially generated. Then $ \mr{FM}(X)=\{ X \}$. 
\end{thm}

\begin{proof}
Let $Y$ be an arbitrary Fourier-Mukai partner of $X$ and $\Phi :D(X) \to D(Y)$ an equivalence. 
We fix $\Phi$. We would like to show that $Y $ is isomorphic to $X$.

Choose a very ample line bundle $L_Y$ on $Y$ and fix it. 
Since the induced morphism $\Phi _* :\Aut (D(X)) \to \Aut (D(Y))$ is an isomorphism, we have 
\[
\exists F \in \Aut (D(X)) \mbox{ s.t. }\Phi _*(F) = L_Y \otimes (-). 
\]
We remark that the following diagram commutes:
\[
\begin{CD}
D(X) @>\Phi>> D(Y)\\
@VFVV @VV L_Y\otimes (-)V \\
D(X)@>\Phi>> D(Y).
\end{CD}
\]
Since $k L_Y$ has a global section for any positive integer $k \in \bb Z_{>0}$, we can make a morphism $E \to E \otimes k L_Y$ for any $E \in D(Y)$. 
Thus, for any $k\in \bb Z_{>0}$ and $E \in D(Y)$, we have $\Hom_{D(Y)}(E,E\otimes k L_Y) \neq 0$. 
Thus we have
\[
\Hom _{D(X)}(\mca O_x, F^{k}(\mca O_x) ) \cong \Hom _{D(Y)}(\Phi (\mca O_x),\Phi (\mca O_x) \otimes kL_Y) \neq 0.
\]

As $\Aut (D(X))$ is trivially generated, $F$ should be written by 
\[
F(?) = L_X \otimes f_* (?) [n],
\]
where $f \in \Aut (X)$, $L_X \in \mr{Pic}(X)$ and $n \in \bb Z$. 
We shall prove that $n=0$ and $f = id _X$. 

Suppose to the contrary that $n \neq 0$. 
Since $n \neq 0$, for sufficiently large $\ell \in \bb Z_{>0}$, we have $\Hom_{D(X)}(\mca O_x,F^{\ell}(\mca O_x)) =0 $ where $F^{\ell}$ is the $\ell $ times composition of $F$.  
This is contradiction. 
Hence $n$ should be $0$. 

We assume that $f \neq id _X$. Then there is a closed point $x \in X$ such that $f(x ) \neq x$. Since $F(\mca O_x) = \mca O_{f(x)}$, we have 
\[
\Hom _{D(X)}(\mca O_x, F(\mca O_x)) =0. 
\]
This is contradiction.  

Thus we have $F = L_X \otimes (-)$. 
Hence for any positive integer $k$, 
\[
\Phi (\mca O_x) \otimes k L_Y = \Phi (\mca O_x \otimes k L_X ) = \Phi (\mca O_x).  
\]
Thus each Hilbert polynomial of $H^i(\Phi (\mca O_x))$ with respect to $L_Y$ is constant. 
Since $L_Y$ is very ample, it follows that $\dim \mr{Supp}(H^i(\Phi (\mca O_x)))=0$. 
Thus $\dim \mr{Supp}(\Phi (\mca O_x))=0$. 
By Lemma \ref{3.1}, we have 
\[
\Phi (\mca O_x) = \mca O_{y_x}[n_x],
\]
for some $y_x \in Y$ and $n_x \in \bb Z$. By the first half assertion of Corollary \ref{2.3}, $n_x$ is locally constant. Hence, $n_x$ is constant. So we put $n_x=n$. 
Then $Y$ is isomorphic to $X$ by the last half assertion of Corollary \ref{2.3}. 
\end{proof}

\begin{rmk}\label{3.3}
The converse of Theorem \ref{3.2} does not hold. 
For instance, there are projective K3 surfaces $X$ with Fourier-Mukai number one (See \cite{HLOY} or \cite{Ogu}). 
On the other hand, 
as is well-knonw by \cite{ST}, the spherical twist $T_S$ by a spherical object\footnote{For example a line bundle on $X$ is a spherical object. } $S \in D(X)$ gives an autoequivalence of $D(X)$ which does not belong to $\mr{Tri}(X)$. 
Thus $\mr{Tri}(X)$ is a proper subgroup of $\Aut (D(X))$. 
\end{rmk}

Let us consider the following three statements for a projective manifold $X$:
\begin{description}
\item[$(A)$] The canonical bundle $K_X$ (or $-K_X$) is ample. 
\item[$(B)$] The autoequivalence group $\Aut (D(X))$ is trivially generated. 
\item[$(C)$] The Fourier-Mukai number of $X$ is one. 
\end{description}
\cite{B-O} proved that $(A) \Rightarrow (B)$ and $(A) \Rightarrow (C)$. 
As we wrote in Remark \ref{3.3}, the converse does not hold. 
Our theorem claims that $(B) \Rightarrow (C)$. 

Now we would like to show that the proposition $(B) \Rightarrow (A)$ does not hold. 

\begin{prop}\label{toda}
Let $X$ be a projective manifold such that $\deg K_X | _C \neq 0$ for any irreducible curve $C \subset X$. 
Then $\Aut (D(X))$ is trivially generated. 
\end{prop}

For instance, let $Y$ be a projective manifold such that $K_Y$ is ample and 
let $X\to Y$ be the blowing up at a point of $Y$. 
Then $X$ satisfies the assumption. 

\begin{proof}
We choose an arbitrary autoequivalence $F\in \Aut (D(X))$ and fix it. 
Since the functor $\otimes K_X [\dim X]$ is the Serre functor, 
the following diagram commutes up to isomorphisms:
\[
\begin{CD}
D(X) @>F>> D(X) \\
@V\otimes K_XVV @VV\otimes K_XV \\
D(X) @>F>> D(X).
\end{CD}
\]
Thus we have 
\begin{equation}
F(\mca O_x) \simeq F(\mca O_x )\otimes K_X \label{kihon}. 
\end{equation}

It suffices to show that $\dim \mr{Supp}(F(\mca O_x))=0$ by Lemma \ref{3.1} and Corollary \ref{2.3}. 
Suppose to the contrary that $\dim \mr{Supp}(F(\mca O_x)) > 0$.  
Then there is an irreducible curve $C$ contained in $\mr{Supp}(F(\mca O_x))$. 
In particular we assume that $C \subset \mr{Supp}(H^i(F(\mca O_x)))$ for some $i \in \bb Z$. 
Now we put $\mca F = H^i(F(\mca O_x))|_C$. 
Notice that $\rank \mca F >0$. 
Thus we have 
\[
\deg \mca F \otimes K_X|_C - \deg \mca F = \rank \mca F \cdot \deg K_X|_C \neq 0
\]
%
On the other hand $\deg \mca F \otimes K_X|_C - \deg \mca F$ should be $0$ by (\ref{kihon}). 
This is contradiction. 
\end{proof}

The next corollary easily follows from Proposition \ref{toda} and Theorem \ref{3.2}. 

\begin{cor}\label{todanocor}
Notations are being as above. 
Then $\mr{FM}(X)=\{ X \}$. 
\end{cor}

\begin{flushleft}
Kotaro Kawatani

Department of Mathematics

Osaka University

Toyonaka 563-0043, Japan

kawatani@cr.math.sci.osaka-u.ac.jp
\end{flushleft}

\end{document}